\documentclass[a4paper,12pt,oneside]{amsart}
\usepackage{geometry}
\usepackage[T1]{fontenc}
\usepackage[french,english]{babel}

\usepackage{multido}
\usepackage{etoolbox}
\usepackage{pdfpages}
\usepackage{amsmath}
\usepackage{amssymb}
\usepackage{amsthm}
\usepackage{mathrsfs}
\usepackage{array}
\usepackage{blkarray}
\usepackage{cancel}
\usepackage{extarrows}
\usepackage{setspace}
\usepackage{graphicx}
\usepackage[all]{xy}
\usepackage{multido}
\usepackage{tabularx}
\usepackage{arydshln}
\usepackage{mathtools}
\usepackage{tikz}
\usetikzlibrary{fit}
\usepackage{tikz-3dplot}
\usepackage{pdflscape}
\usepackage{multirow}
\usetikzlibrary{calc,backgrounds,matrix}
\usetikzlibrary{arrows,shapes,patterns}
\usepackage{supertabular}
\usepackage{titletoc}
\usepackage{fancybox}
\usepackage{fancyhdr}
\usepackage{enumitem}
\usepackage{diagbox}
\usepackage{nicematrix}
\usepackage{hyperref}

\setenumerate[1]{itemsep=-3pt, label=$(\mathrm{\arabic*})$}

\newcommand{\cF}{\mathcal{F}}

\newtheorem{Def}{\bf Definition}
\newtheorem{Thm}{\bf Theorem}
\newtheorem{ThmIntro}{\bf Theorem}
\newtheorem{Prop}{\bf Proposition}

\newtheorem{Rem}{\bf Remark}
\newtheorem{Lem}{\bf Lemma}

\edef\myindent{\the\parindent}

\setlist[itemize,1]{label=\textcolor{black}{\small $\bullet$}}
\setlist[itemize,2]{label=\textcolor{black}{\tiny $\bullet$}}
\setlist[itemize,3]{label=$\diamond$}
\setlist[itemize,4]{label=$\bullet$}

\begin{document}
\title{Construction of Transverse Metrics from Transversely Elliptic Operators on Foliated Manifolds}
\author{Wenran Liu\\Civil Aviation University of China}
\address{Wenran Liu, Sino-European Institute of Aviation Engineering, Civil Aviation University of China, 2898 Jinbei Street Dongli District, 300300 Tianjin, P.R.China }
\email{wrliu@cauc.edu.cn}
\thanks{This research project is supported by the Starting Research Fund from Civil Aviation University of China, Project No. 2020KYQD58
}
\date{}

\maketitle
\begin{abstract}
The existence of a transverse metric on a foliated manifold is typically used as a framework for studying transversely elliptic basic differential operators. This article details a method for constructing such a transverse metric directly from the principal symbol of the operator. Consequently, rather than treating the transverse metric as a mere working assumption, we demonstrate that its existence is a necessary condition for the very presence of a transversely elliptic operator.
\end{abstract}

\selectlanguage{english}
\tableofcontents

\normalsize
\section{Introduction}
The notion of bundle-like metric on foliated manifold was first introduced by Reinhart \cite{Rei}. A foliated manifold equipped with a bundle-like metric is called a Riemannian foliation. An equivalent definition is the existence of a transverse (holonomy-invariant) metric on the normal bundle of the foliation. In \cite{Mol, Mol0}, Molino proposed a theory for the study of Riemannian foliation. In \cite{Bre} and in \cite{Car}, Carrière and Breuillard-Gelander studied the growth of Riemannian foliations.  In \cite{Ka}, El.Kacimi-Alaoui introduced the notion of transversely elliptic operator and proved the association of a transversely elliptic operator over a Riemannian foliation and an equivariant operator over a compact manifold constructed by Molino's theory. In \cite{BKR}, J.Bruning, F.Kamber and K.Richardson proposed an index theory for basic Dirac operators on Riemannian foliations. In \cite{GL}, A.Gorokhovsky and J.Lott gave an index formula for a transverse Dirac-type operator in case of Abelian Molino sheaf. As their work are all based on the Riemannian assumption of the foliation, El.Kacimi-Alaoui asked in \cite[Question2.9.2]{AKHR} whether a transversely elliptic operator exists only if the foliation is Riemannian.\smallskip\\
\indent In the present work,we answer this question in a more general case. Let us recall El.Kacimi-Alaoui's question as follows:
\smallskip
\begin{center}
\begin{minipage}{0.8\textwidth}
Let $M$ be a compact manifold with a foliation $\cF$ which admits non-constant basic functions. Suppose that there exists a transversely elliptic operator acting on these functions. Is the foliation $\cF$ Riemannian?
\end{minipage}
\end{center}
\smallskip
\indent In \cite{Liu}, an averaging method is applied to construct a basic connection on a foliated bundle over a Riemannian foliation. We consider the holonomy groupoid restricted to a complete transversal, whose closure is shown to be a Hausdorff, proper Lie groupoid equipped with a Haar system. By combining this structure with suitable cut-off functions, the basic connection is successfully constructed by averaging an initial foliated connection. \smallskip\\
\indent Unfortunately, in the absence of a Riemannian foliation structure, the average method is no longer applicable. Even the holonomy pseudogroup approach \cite{Hae,Sac,Sal,Sal1} introduced by Haefliger fundamentally relies on the foliation being Riemannian.
\\
\indent Furutani \cite{Fu} investigated the transformation group of a compact manifold that commutes with an elliptic differential operator. Utilizing global analysis on manifolds, his main result establishes that, when equipped with the compact-open topology, the group of such commuting diffeomorphisms forms a compact Lie group. However, extending this work to the foliated setting is inherently difficult without first ensuring the existence of a transverse metric. While a Riemannian metric always exists on any smooth manifold, a transverse metric does not necessarily exist, a classic counterexample being the Reeb foliation \cite{CC}, another counterexample is constructed in \cite{Liu}. Indeed, a transverse metric is a global geometric object, and its holonomy-invariant property is far more restrictive than a mere collection of local holonomy-invariance conditions within individual foliated charts.
\smallskip\\
The answer to the question raised by El Kacimi-Alaoui is positive. In this paper, we construct a transverse metric directly from the principal symbol of a transversely elliptic operator, using the symmetric tensor algebra as our primary tool. Given a transversely elliptic basic differential operator, we first consider its composition with its conjugate. This composition remains transversely elliptic, with the key advantage of having real-valued coefficients and an even order. Its principal symbol can then be viewed as a positive-definite section of the symmetric $(2m,0)$-tensor bundle over the normal bundle of the foliated manifold. Our main method consists in constructing a positive-definition section of the symmetric $(0,2)$-tensor bundle over the normal bundle, and proving its holonomy invariance.Let us list out main result below. 
\begin{ThmIntro}
Let $M$ be a foliated manifold with a foliation $\cF$ which admits non-constant basic functions. If there exists a transversely elliptic basic differential operator acting on the space of basic functions, then $\cF$ is a Riemannian foliation. 
\end{ThmIntro}
\begin{ThmIntro}
Each second order transversely elliptic basic differential operator acting on the space of basic functions with real-valued basic function coefficients is a generalized transverse Laplace operator. 
\end{ThmIntro}
\begin{ThmIntro}
Let $M$ be a foliated manifold with a foliation $\cF$ which admits non-constant basic functions, $\mathcal{E}\rightarrow M$ be an $\cF$-bundle. If there exists a transversely elliptic basic differential operator on $\mathcal{E}$, then $\cF$ is a Riemannian foliation.
\end{ThmIntro}

The paper is organized as follows. In Section 2, we study the symmetric tensor algebra of a vector space and detail a method for constructing a symmetric $(0,2)$-tensor from a given positive-definite symmetric $(2m,0)$-tensor. In Section 3, we extend this algebraic construction to symmetric tensor bundles over manifolds. Section 4 is devoted to foliated manifolds equipped with a transversely elliptic basic differential operator. More precisely, in Section 4.1, we establish the existence of a transverse metric for the case where the transversely elliptic operator acts on the space of basic functions; notably, we show that the compactness assumption in El Kacimi-Alaoui's question can be removed. Finally, in Section 4.2, we generalize this result to a broader setting: we prove that if a transversely elliptic basic differential operator is defined on an $\mathcal{F}$-bundle, then the foliation must necessarily be Riemannian.\\[0.5em]

\noindent {\bf\em Acknowledgements:} \\
\noindent The author is grateful to Moulay Tahar Benameur, Paul-Emile Paradan and Alexandre Baldare for various helpful discussions on this project. In particular, the author deeply acknowledges Professor Ma Xiaonan for his encouragement and continuous support. This work is dedicated to the memory of Professor Moulay Tahar Benameur, my deeply respected doctoral co-advisor. I am profoundly saddened by his passing. Our last meeting in 2020 at the Université de Montpellier will remain forever engraved in my memory.\\
\noindent The author welcomes and gratefully accepts all remarks, comments, doubts and criticisms regarding this paper. \\


\section{Symmetric tensor space and construction method}
Throughout this artical, let $m,n\in\mathbb{N}^*$. Let $S^{2m,0}(\mathbb{R}^n)$ be the space of $(2m,0)$-type symmetric tensors on $\mathbb{R}^n$. An element of $S^{2m,0}(\mathbb{R}^n)$ is a real coefficient homogeneous polynomial on $\mathbb{R}^n$ of degree $2m$. Each $Q\in S^{2m,0}(\mathbb{R}^n)$ corresponds to a function $\widetilde{Q}$ from $\mathbb{R}^n$ to $\mathbb{R}$
$$
\widetilde{Q}(X)=\sum_{|\alpha|=2m}a_{\alpha} x_1^{\alpha_1}x_2^{\alpha_2}\cdots x_n^{\alpha_n}
$$
where $\alpha=(\alpha_1,\alpha_2,\cdots,\alpha_n)\in \mathbb{N}^n$, $|\alpha|=\alpha_1+\alpha_2+\cdots+\alpha_n$, $a_\alpha\in\mathbb{R}$ and the variable $X=[x_1,x_2,\cdots,x_n]^T$. Here, we identify the vector space $\mathbb{R}^n$ with real vertical matrices of order $n$ and $T$ means the matrix transposition.
\begin{Def}
We say that $Q\in S^{2m,0}(\mathbb{R}^n)$ is
\begin{enumerate}
\item[$\bullet$] positive if $\forall\,X\in\mathbb{R}^n$, $\widetilde{Q}(X)\geq 0$ 
\item[$\bullet$] definite-positive if it is positive and $\widetilde{Q}(X)=0$ if and only if $X=0$.
\end{enumerate}
\end{Def}
Denote by $S^{2m,0}_{+}(\mathbb{R}^n)$ (resp. $S^{2m,0}_{++}(\mathbb{R}^n)$) the set of positive (resp. definite-positive) real coefficient homogeneous polynomial on $\mathbb{R}^n$ of degree $2m$. Remark that when $n=1$, the positiveness and the definite positiveness coincide. The following proposition is clear.
\begin{Prop}\label{Prop2.1}
For each matrix $A\in \mathrm{GL}_n(\mathbb{R})$, the right mulplication by $A$:
$$
\mathrm{R}_A:S^{2m,0}(\mathbb{R}^n)\rightarrow S^{2m,0}(\mathbb{R}^n), Q\mapsto \mathrm{R}_AQ
$$
is a bijection where $(\mathrm{R}_AQ)(X)=Q(AX)$. Restricting to $S^{2m,0}_{+}(\mathbb{R}^n)$ and to $S^{2m,0}_{++}(\mathbb{R}^n)$, $R_A$ is also a bijection respectively.
\end{Prop}


By convention, denote by $\mathbb{B}_n^*$ the decentered open $n$-ball, by $\mathbb{S}^{n-1}$ the $(n-1)$-sphere. 
\begin{Prop}
For each $Q\in S^{2m,0}_{++}(\mathbb{R}^n)$,
\begin{enumerate}
\item[$\bullet$] $\{0<Q<1\}=\{X\in\mathbb{R}^n\mid\,0<\widetilde{Q}(X)<1\}$
is a smooth manifold diffeomorphic to $\mathbb{B}_n^*$;
\item[$\bullet$] $\{Q=1\}=\{X\in\mathbb{R}^n\mid\,\widetilde{Q}(X)=1\}$ is a smooth manifold diffeomorphic to $\mathbb{S}^{n-1}$.
\end{enumerate}
\end{Prop}
\begin{proof}
By the compactness of $\mathbb{S}^{n-1}$, the maximum and the minimum of $\widetilde{Q}$ on $\mathbb{S}^{n-1}$ is attained. Since $\widetilde{Q}$ takes positive value on $\mathbb{S}^{n-1}$,
$$
\max_{\mathbb{S}^{n-1}}\widetilde{Q}\geq \min_{\mathbb{S}^{n-1}}\widetilde{Q}>0.
$$
The map $\theta\mapsto \left[\widetilde{Q}(\theta/\|\theta\|)\right]^{-\frac{1}{2m}} \theta$ is a well-defined diffeomorphism from $\mathbb{B}_n^*$ to $\{0<Q<1\}$, is also a well-defined diffeomorphism from $\mathbb{S}^{n-1}$ to $\{Q=1\}$, where $\|\cdot\|$ is the canonical Euclidean norm.
\end{proof}
Let $d\lambda$ be the Lebesgue measure on $\mathbb{R}^n$. 
\begin{Lem}\label{Lem2.1}
For each $Q\in S^{2m,0}_{++}(\mathbb{R}^n)$, 
$$
\int_{b\in \{0<Q<1\}} d\lambda >0.
$$
\end{Lem}
\begin{proof}
By change of variables $b=\left[\widetilde{Q}(\theta/\|\theta\|)\right]^{-\frac{1}{2m}}\theta$, $\theta\in\mathbb{B}_n^*$,
$$
\int_{b\in \{0<Q<1\}}d\lambda\vert_b= \int_{\theta\in\mathbb{B}_n^*}\left[\widetilde{Q}(\theta/\|\theta\|)\right]^{-\frac{n}{2m}}d\lambda\vert_\theta\geq \left(\max_{\mathbb{S}^{n-1}}\widetilde{Q}\right)^{-\frac{n}{2m}}\frac{\pi^{\frac{n}{2}}}{\Gamma(\frac{n}{2}+1)}>0,
$$
where $\Gamma$ is Euler's gamma function.
\end{proof}
\begin{Lem}\label{Lem2.2}
$S^{2m,0}_{++}(\mathbb{R}^n)$ is a smooth manifold of dimension $\frac{(2m+n-1)!}{(2m)!(n-1)!}$.
\end{Lem}
\begin{proof}
Each element of $S^{2m,0}(\mathbb{R}^n)$ is completely and uniquely determined by its coefficient $(a_\alpha)_{\alpha\in\mathbb{N}^n,|\alpha|=2m}$, so $S^{2m,0}(\mathbb{R}^n)$ has a natural smooth manifold structure diffeomorphic to $\mathbb{R}^{\frac{(2m+n-1)!}{(2m)!(n-1)!}}$. It suffices to show that $S^{2m,0}_{++}(\mathbb{R}^n)$ is an open subset of $S^{2m,0}(\mathbb{R}^n)$. \\
Indeed, by homogeneity, $Q\in S^{2m,0}_{++}(\mathbb{R}^n)$ if and only if $\min\limits_{\mathbb{S}^{n-1}} \widetilde{Q}>0$. Then, for each $Q\in S^{2m,0}_{++}(\mathbb{R}^n)$ with $ \widetilde{Q}=\sum_{|\alpha|=2m}a_\alpha {x_1}^{\alpha_1}{x_2}^{\alpha_2}\cdots {x_n}^{\alpha_n}$, there exists $\delta>0$ such that $\forall\, (a^\prime_\alpha)_{\alpha\in\mathbb{N}^n,|\alpha|=2m}$ satisfying $\|a^\prime_\alpha-a_\alpha\|<\delta$, we have
$$
\min_{S^{n-1}} \Big(\sum_{|\alpha|=2m}a^\prime_\alpha {x_1}^{\alpha_1}{x_2}^{\alpha_2}\cdots {x_n}^{\alpha_n}\Big)>\frac{1}{2} \min_{S^{n-1}} Q>0.
$$
\end{proof}
Denote by $\mathcal{S}_n^{+}$ (resp. $\mathcal{S}_n^{++}$) the set of positive (resp. definite-positive) symmetric matrix of order $n$. 
\begin{Prop}\label{Prop2.4}
For each $Q\in S^{2m,0}_{++}(\mathbb{R}^n)$ and $b\in \{0<Q<1\}$, the equation
\begin{equation}\label{eq:2.4}
\mathbb{E}_{Q,b}:\widetilde{Q}(X)=1\ \text{with}\ X \parallel b
\end{equation}
admits a unique solution $C_{Q,b}$ up to a minus sign. The matrix $C_{Q,b}{C_{Q,b}}^T$ is well-defined and belongs to $\mathcal{S}_n^+$ whose kernel is the subspace $b^\perp=\{X\in\mathbb{R}^n\mid\, b^T X=0\}$ of $\mathbb{R}^n$.
\end{Prop}
\begin{proof}
The solution is $C_{Q,b}=\pm  \left[\widetilde{Q}(b)\right]^{-\frac{1}{2m}} b$. Indeed, since $Q$ is homogeneous of degree $2m$, 
$$
\widetilde{Q}(C_{Q,b})=\left[\widetilde{Q}(b)^{-1/2m}\right]^{2m} \widetilde{Q}(b)=1.
$$
Up to a minus sign, $C_{Q,b}$ is unique. Since the minus sign cancels, the symmetric matrix $C_{Q,b}{C_{Q,b}}^T$ is well-defined and positive. we calculate its kernel
$$
\mathrm{Ker}(C_{Q,b}{C_{Q,b}}^T)=\left\{X\in\mathbb{R}^n\big|\, {C_{Q,b}}^TX=0\right\}=\left\{X\in\mathbb{R}^n\big|\, b^T X=0\right\}.
$$
\end{proof}
We define the map
$$
\begin{array}{rcl}
\phi:S^{2m,0}_{++}(\mathbb{R}^n)&\rightarrow&\mathcal{S}_n^{++}\\
Q&\mapsto&\phi(Q)=n\left[\displaystyle\int_{b\in \{0<Q<1\}} d\lambda\right]^{-1}\displaystyle\int_{b\in \{0<Q<1\}} C_{Q,b}C_{Q,b}^T d\lambda
\end{array}
$$

\begin{Prop}\label{Prop2.4}
The map $\phi$ is well-defined and smooth. For each $Q\in S^{2m,0}_{++}(\mathbb{R}^n)$, $A\in\mathrm{GL}_n(\mathbb{R})$, we have 
$$
\phi(Q)=A\,\phi(\mathrm{R}_AQ)\,A^{T}.
$$
where $\mathrm{R}_A$ is the right multiplication by $A$ in Proposition \ref{Prop2.1}
\end{Prop}
\begin{proof}
By Lemma \ref{Lem2.1}, the denominator of $\phi(Q)$ is positive. Let us prove the definite-positiveness of the numerator of $\phi(Q)$. Suppose not, $\exists\, X_0\in\mathbb{R}^n$ with $X_0\neq 0$ such that
$$
\int_{b\in \{0<Q<1\}} X_0^{T}C_{Q,b}{C_{Q,b}}^TX_0 d\lambda=0.
$$
By the positiveness of $C_{Q,b}{C_{Q,b}}^T$, for every $b\in\{0<Q<1\}$, 
$$
X_0^{T}C_{Q,b}{C_{Q,b}}^TX_0=0.
$$
This implies $X_0\in b^\perp$ for any $b$. In particular, take $b$ to be $b_0 = 2^{\frac{1}{2m}} \left[\widetilde{Q}(X_0)\right]^{-\frac{1}{2m}}X_0$. It is easy to check that $b_0\in\{0<Q<1\}$. But we have a contradicition that
$$
b_0^TX_0=[Q(X_0)/2]^{-1/2m} X_0^TX_0>0.
$$
It is known that $\mathcal{S}_n^{++}$ is a connected smooth manifold. By lemma \ref{Lem2.2}, $S^{2m,0}_{++}(\mathbb{R}^n)$ is a smooth manifold. The smoothness of $\phi$ comes from that of the set $\{0<Q<1\}$ and the solution $C_{Q,b}$ for any fixed $b$ with respect to the variable $Q$, or equivalently the coefficients of $Q$.\\
For each $Q\in S^{2m,0}_{++}(\mathbb{R}^n)$, $A\in\mathrm{GL}_n(\mathbb{R})$, $X\mapsto A^{-1}X$ is a diffeomorphism from $\{0<Q<1\}$ to $\{0<\mathrm{R}_AQ<1\}$, meanwhile, is also a diffeomorphism from $\{Q=1\}$ to $\{\mathrm{R}_AQ=1\}$. Then, for every $b\in\{0<Q<1\}$, the equation $\mathbb{E}_{R_AQ,A^{-1}b}$ makes sense. We claim that $A^{-1}C_{Q,b}$ solves the equation $\mathbb{E}_{R_AQ,A^{-1}b}$. Indeed, since $C_{Q,b}$ solves $\mathbb{E}_{Q,b}$,
$$
\widetilde{\mathrm{R}_AQ}(A^{-1}C_{Q,b})=\widetilde{Q}\big(A A^{-1}C_{Q,b}\big)=1\ \text{and}\ (A^{-1}C_{Q,b}) \parallel (A^{-1}b).
$$
By the uniqueness of solution up to a minus sign, we have
$$
C_{Q,b}=A C_{\mathrm{R}_AQ,A^{-1}b}.
$$
Take change of variables $b^\prime=A^{-1}b$, we have
$$
\begin{array}{rl}
\phi(Q)=&n\left[\displaystyle\int_{b\in \{0<Q<1\}} d\lambda\vert_b\right]^{-1}\displaystyle\int_{b\in \{0<Q<1\}} C_{Q,b}{C_{Q,b}}^T d\lambda\vert_b\\[1em]
=&n\left[\displaystyle\int_{b^\prime\in \{0<\mathrm{R}_AQ<1\}} \cancel{(\det A)} d\lambda\vert_{b^\prime}\right]^{-1}\displaystyle\int_{b^\prime\in \{0<\mathrm{R}_AQ<1\}} A C_{\mathrm{R}_AQ,b^\prime}{C_{\mathrm{R}_AQ,b^\prime}}^T A^T \cancel{(\det A)} d\lambda\vert_{b^\prime}\\[1em]
=&A\phi(\mathrm{R}_AQ)A^T.
\end{array}
$$
\end{proof}
Let us see the special case of $m=1$.
\begin{Prop}\label{Prop2.5}
For each $[a^{ij}]\in \mathcal{S}_n^{++}$, if $Q\in S^{2,0}(\mathbb{R}^n)$ is defined by $\widetilde{Q}(X)=\sum_{i,j=1}^n a^{ij}x_ix_j$, then $Q\in S^{2,0}_{++}(\mathbb{R}^n)$ and
$$
\phi(Q)=[a^{ij}]^{-1}.
$$
\end{Prop}
\begin{proof}
The definite-positiveness of $Q$ is equivalent to that of $[a^{ij}]$. For the second assertion, by Proposition \ref{Prop2.4}, it suffices to prove that if $[a^{ij}]=I_n$, equivalently $\widetilde{Q}(X)=\sum_{k=1}^n x_k^2$, then $\phi(Q)=I_n$. We leave the cases of $n=1$ or $2$ to readers. When $n\geq 3$, take the hyperspherical coordinates $(\rho,\theta)\in (0,+\infty)\times \mathbb{S}^{n-1}$ where $\theta=(\theta_1,\cdots,\theta_{n-1})\in [0,\pi]^{n-2}\times [0,2\pi]$, then $d\lambda = \rho^{n-1}d\rho dS$ where $dS=\sin^{n-2}(\theta_1)\cdots\sin(\theta_{n-2})d\theta_1\cdots d\theta_{n-1}$ is the surface area form of $\mathbb{S}^{n-1}$. It is easy to get 
$$
{C_{Q,\rho\theta}}^T=\left[\cos\theta_1,\sin\theta_1\cos\theta_2,\cdots,\sin\theta_1\cdots\sin\theta_{n-2}\cos\theta_{n-1},\sin\theta_1\cdots\sin\theta_{n-2}\sin\theta_{n-1}\right].
$$
We calculate
$$
\begin{array}{rl}
\phi(Q) =&n\left[\displaystyle\int_{\theta\in \mathbb{S}^{n-1}}\int_{0}^{1}\rho^{n-1}d\rho dS\right]^{-1} \displaystyle\int_{\theta\in \mathbb{S}^{n-1}}\int_{0}^{1} C_{Q,\rho\theta}{C_{Q,\rho\theta}}^T \rho^{n-1}d\rho dS\\ \\
=&n\left[\displaystyle\int_{\theta\in \mathbb{S}^{n-1}}dS\right]^{-1} \displaystyle\int_{\theta\in \mathbb{S}^{n-1}}C_{Q,\rho\theta}{C_{Q,\rho\theta}}^T dS=I_n
\end{array}
$$
The last step is obtained by tedious calculation of each component.
\end{proof}

Next, we pass to the construction on the symmetric tensor space of vector space. Let $V$ be a $n$-dimensional real vector space and $V^*$ be its dual space. Consider $S^{2m,0}(V)$ the space of $(2m,0)$-type symmetric tensors on $V$. Each $P\in S^{2m,0}(V)$ corresponds to a function $\widetilde{P}:V^*\rightarrow\mathbb{R}$. Precisely, chosen a basis $\mathscr{B}=(e_1,e_2,\cdots,e_n)$ of $V$, $P$ has the expression
$$
P = \sum_{|\alpha|=2m} a_{\alpha}e_1^{\alpha_1}e_2^{\alpha_2} \cdots e_n^{\alpha_n},
$$
where $\alpha=(\alpha_1,\cdots,\alpha_n)\in\mathbb{N}^n$, $|\alpha|=\alpha_1+\cdots+\alpha_n$ and $a_\alpha\in\mathbb{R}$, then for any covector $v^*\in V^*$, 
$$
\widetilde{P}(v^*) = \sum_{|\alpha|=2m} a_{\alpha}\big<e_1,v^*\big>^{\alpha_1}\big<e_2,v^*\big>^{\alpha_2} \cdots \big<e_n,v^*\big>^{\alpha_n}
$$
where $\big<\cdot,\cdot\big>$ denotes the contraction of vector and covector. 
\begin{Def}
We say that $P\in S^{2m,0}(V)$ is 
\begin{enumerate}
\item[$\bullet$] positive if $\forall\,v^*\in V^*$, $\widetilde{P}(v^*)\geq 0$,
\item[$\bullet$] definite-positive if $P$ is positive and $\widetilde{P}(v^*)=0$ if and only if $v^*=0$. 
\end{enumerate}
Denote by $S^{2m,0}_{+}(V)$ (resp. $S^{2m,0}_{++}(V)$) the set of positive (resp. definite-positive) symmetric $(2m,0)$-type tensors on $V$.
\end{Def}
One may associate $P$ with $P_{\mathscr{B}}\in S^{2m,0}(\mathbb{R}^n)$ under a chosen basis $\mathcal{B}$ as above by
$$
\widetilde{P_{\mathscr{B}}}(X) = \widetilde{P}(x_1e_1^*+\cdots+x_ne_n^*)= \sum_{|\alpha|=2m} a_{\alpha}x_1^{\alpha_1}x_2^{\alpha_2} \cdots x_n^{\alpha_n}.
$$ 
where $\mathscr{B}^*=(e_1^*,e_2^*,\cdots,e_n^*)$ is the dual basis of $\mathscr{B}$. From now on, we denote $\widetilde{P_\mathscr{B}}$ simply by $P_\mathscr{B}$. Evidently, $P\in S^{2m,0}_{+}(V)$ (resp, $S^{2m,0}_{++}(V)$) if and only if $P_{\mathscr{B}}\in S^{2m,0}_{+}(\mathbb{R}^n)$ (resp, $S^{2m,0}_{++}(\mathbb{R}^n)$).\\
Analoguely, we have $S^{0,2}_{++}(V)$ the set of definite-positive symmetric $(0,2)$-type tensors on $V$. Clearly, each element of $S^{0,2}_{++}(V)$ is an inner product on $V$. 
\begin{Def}
Define the map
$$
\begin{array}{rcl}
\Phi_V : S^{2m,0}_{++}(V) &\rightarrow &S^{0,2}_{++}(V)\\
P&\mapsto&\sum_{i,j=1}^n \big[\phi(P_{\mathscr{B}})\big]_{ij} e_i^* e_j^*
\end{array}
$$
where $\big[\phi(P_{\mathscr{B}})\big]_{ij}$ is the $(i,j)$-component of the matrix $\phi(P_{\mathscr{B}})$.
\end{Def}
\begin{Prop}\label{Prop2.6}
The map $\Phi_V$ is well-defined. 
\end{Prop}
\begin{proof}
It suffices to show that $\Phi_V(P)$ does not depend on the choice of basis of $V$. Let $\mathscr{C}=(\varepsilon_1,\varepsilon_2,\cdots,\varepsilon_n)$ be another basis of $V$ and $\mathscr{C}^*=(\varepsilon_1^*,\varepsilon_2^*,\cdots,\varepsilon_n^*)$ be the dual basis. Denote by $A$ the change matrix from $\mathscr{C}$ to $\mathscr{B}$, i.e.
$e_l=\sum_{k=1}^n A_{kl}  \varepsilon_k$. It only needs to prove
$$
\phi_{\mathscr{B}}(P)=A^T \phi_{\mathscr{C}}(P) A.
$$
As covector on $V^*$ with coordinates $X=[x_1,\cdots,x_n]^T$ under $\mathscr{C}^*$ corresponds to that with coordinates $A^TX$ under $\mathscr{B}^*$, so we have
\begin{equation}\label{eq:2.1}
P_{\mathscr{C}}=\mathrm{R}_{A^T}P_{\mathscr{B}}.
\end{equation}
Apply Proposition \ref{Prop2.4}, the proof is complete.
\end{proof}
Let us see the specital case of $m=1$.
\begin{Prop}
For each $[g^{ij}]\in \mathcal{S}_n^{++}$, if $P\in S^{2,0}(V)$ is defined by $\widetilde{P}=\sum_{i,j=1}^n g^{ij}e_i e_j$ under some basis $\mathscr{B}=(e_1,\cdots,e_n)$, then $P\in S^{2,0}_{++}(V)$ and 
$$
\Phi_V(P) = \sum_{i,j=1}^n g_{ij} e_i^*e_j^*
$$
where $[g_{ij}]$ is the inverse matrix of $[g^{ij}]$ and $\mathscr{B}^*=(e_1^*,\cdots,e_n^*)$ is the dual basis of $\mathscr{B}$.
\end{Prop}
\begin{proof}
Apply Proposition \ref{Prop2.5} to $P_{\mathscr{B}}$.
\end{proof}
\section{Construction on symmetric tensor bundles}
Let $M$ be a smooth manifold and $E$ be a real vector bundle over $M$ of rank $n$. Consider $S^{2m,0}(E)$, $S^{0,2}(E)$ the symmetric $(2m,0)$-type, $(0,2)$-type tensor bundle on $M$, respectively. For tensor bundles, see \cite[Section5.2]{AMR}. Let
$$
S^{2m,0}_{++}(E)=\bigsqcup_{z\in M}S^{2m,0}_{++}(E\vert_z),\ S^{0,2}_{++}(E)=\bigsqcup_{z\in M} S^{0,2}_{++}(E\vert_z)
$$
%
\begin{Lem}
We have
\begin{enumerate}
\item[$\bullet$] $S^{2m,0}_{++}(E)$ is a fiber bundle on $M$ whose fiber is diffeomorphic to $S^{2m,0}_{++}(\mathbb{R}^n)$; 
\item[$\bullet$] $S^{0,2}_{++}(E)$ is a fiber bundle on $M$ whose fiber is diffeomorphic to $S^{0,2}_{++}(\mathbb{R}^n)$.
\end{enumerate}
\end{Lem}
\begin{proof}
On each local chart $U$, we have a local trivilization $E\vert_U\simeq U\times \mathbb{R}^n$, which induces a lcoal trivilization $S^{2m,0}(E)\vert_U\simeq U\times S^{2m,0}(\mathbb{R}^n)$. By Lemma \ref{Lem2.2}, $S^{2m,0}_{++}(\mathbb{R}^n)$ has a manifold structure which is preserved under $\mathrm{GL}_n(\mathbb{R})$, therefore $S^{2m,0}_{++}(E)$ is a fiber bundle over $M$ whose fiber is diffeomorphic to $S^{2m,0}_{++}(\mathbb{R}^n)$. Analogue argument is true for $S^{0,2}_{++}(E)$.
\end{proof}
\begin{Prop}\label{Prop3.1}
The following map
$$
\begin{array}{rcl}
\Psi_E:C^\infty\left(M,S_{++}^{2m,0}(E)\right)&\rightarrow&C^\infty\left(M,S_{++}^{0,2}(E)\right)\\
P&\mapsto&\Psi_E(P)
\end{array}
$$
is well-defined bundle morphism where $\Psi_E(P)\vert_z=\Phi_{E\vert_z}(P\vert_z)$.
\end{Prop}
\begin{proof}
It only needs to prove locally, for each local chart $U$ where $E\vert_U\simeq U\times \mathbb{R}^n$, let $\mathscr{B}=(e_1,e_2,\cdots,e_n)$ be the smooth local frame of $E\vert_U$ induced by the canonical basis of $\mathbb{R}^n$. Denote by $\mathscr{B}^*=(e_1^*,e_2^*,\cdots,e_n^*)$ the dual frame. For each $z\in U$, we have
$$
\Psi_E(P)\vert_z=\Phi_{E\vert_z}(P\vert_z)= \sum_{i,j=1}^n \left[\phi\Big((P\vert_z)_{\mathscr{B}}\Big)\right]_{ij}e_i^*e_j^*.
$$
By the smoothness of the dual frame $\mathscr{B}^*$ and that of $\phi$ proved in Proposition \ref{Prop2.4}, it suffices to prove the smoothness of $(P\vert_z)_{\mathscr{B}}\in S^{2m,0}_{++}(\mathbb{R}^n)$ with respect to $z$. This is true because the section $P$ and the local frame $\mathscr{B}$ are smooth. 
\end{proof}

\section{Construction of transverse metric on foliated manifolds}
Let $M$ be a smooth manifold of dimension $n$, $\cF$ be a foliation of codimension $q$ on $M$. We may assume $M$ to be connected, or we do constuction on its each connected component. The foliation $\cF$ can be regarded as a integrable distribution of subbundle of rank $(n-q)$ of the tangent bundle $TM$. We say simply that $(M,\cF)$ is a foliated manifold. Recall some basic notions of foliation theory. For detail, see \cite{CC,Hec,Mol,Moe}. The normal bundle $\nu \cF$ is the quotient bundle $TM\slash \cF$. Denote by $\frak{X}(M)$ the space of vector fields on $M$, by $\frak{X}(\cF)$ the Lie algebra of vector fields tangent to $\cF$. We say that $Y\in\frak{X}(M)$ is a foliated vector field on $M$ if $\forall\, X\in\frak{X}(F)$, $\mathcal{L}_X Y=0$, where $\mathcal{L}$ is the Lie derivative. Denote by $\frak{X}(M,\cF)$ the space of foliated vector fields on $M$, quotient by $\frak{X}(\cF)$ to get the space of basic fields $l(M,\cF)=\frak{X}(M,\cF)\slash\frak{X}(\cF)$. All these notions above can be defined locally on a foliated chart.
\begin{Def}
A transverse metric on the foliated manifold $(M,\cF)$ is a positive $C^\infty(M)$-bilinear form $g$ on $\frak{X}(M)$ such that
\begin{enumerate}
\item[(i)] $\mathrm{Ker}(g\vert_z)=\cF\vert_z$, for any $z\in M$
\item[(ii)] $\mathcal{L}_X g$ = 0 for any $X\in\frak{X}(\cF)$.
\end{enumerate}
A foliation $\cF$ of $M$ together with a transverse metric $g$ on $(M,\cF)$ is called a Riemannian foliation on $M$.
\end{Def}
\begin{Rem}
The second property in the defintion of transverse metric is called the holonomy invariance. Each transverse metric $g$ on $(M,\cF)$ induces a definite-positive $C^\infty(M)$-bilinear form $g^{\nu\cF}$ on $C^\infty(M,\nu\cF)$ satisfying the holonomy invariance $\mathcal{L}_X g^{\nu\cF} = 0$ for any $X\in\frak{X}(\cF)$. Conversely, if we have constructed such a form, a transverse metric on $(M,\cF)$ can be easily obtained by the pullback map $(\nu\cF)^*\rightarrow T^*M$.
\end{Rem}
\subsection{Transversely Elliptic operators acting on basic functions}
A complex-valued smooth function $f\in C^\infty(M)$ on $M$ is said to be basic if $\forall X\in \frak{X}(\cF)$, $X(f)=0$. Denote the space of complex-valued basic function on $M$ by $C_b(M)$. Assume that the foliation $\cF$ admits non-constant basic functions, equivalently $C_b(M)$ is not isomorphic to $\mathbb{C}$.  Since $\cF$ is a subbundle of $TM$, we have a natural isomorphism between the transverse cotangent bundle $T^*M\slash \cF^*$ and the dual bundle $(\nu\cF)^*$ of the normal bundle $\nu\cF$.

\begin{Def}
A basic differential operator $L$ of order $m\in\mathbb{N}^*$ acting on $C_b(M)$ is a map $L:C_b(M)\rightarrow C_b(M)$ such that in each local foliated chart of coordinates $(x,y)=(x_1,\cdots,x_{n-q},y_1,\cdots,y_q)$, $L$ has the expression
$$
L=\sum_{|\alpha|\leq m} a_{\alpha}(y) \dfrac{\partial^{|\alpha|}}{\partial {y_1}^{\alpha_1} \partial {y_2}^{\alpha_2}\cdots \partial {y_q}^{\alpha_q}}
$$ 
where $\alpha=(\alpha_1,\cdots,\alpha_q)\in \mathbb{N}^q$, $|\alpha|=\alpha_1+\cdots+\alpha_q$ and all coefficients $a_{\alpha}$ are complex-valued basic functions. The principal symbol of $L$ at the point $z$ and the covector $\xi\in T^*_z M$ is 
$$
\sigma(L)(z,\xi)=\sum_{|\alpha|=m} a_{\alpha}(y) {\xi_1}^{\alpha_1}\cdots {\xi_q}^{\alpha_q}.
$$
We say that $L$ is transversely elliptic if $\sigma(L)(z,\xi)$ is non-zero for every $z\in M$ and every transverse covector $\xi\in (T^*M\slash \cF^*)\vert_z$ different from $0$.
\end{Def}

\begin{Thm}
Let $M$ be a foliated manifold with a foliation $\cF$ which admits non-constant basic functions. If there exists a transversely elliptic basic differential operator acting on the space of basic functions, then $\cF$ is a Riemannian foliation. 
\end{Thm}
\begin{proof}
Let $L$ be such an operator of order $m\in\mathbb{N}^*$. Consider $L\circ \overline{L}$ where $\overline{L}$ is the conjugate operator (taking conjugate of coefficients of $L$). Evidently, $L\circ \overline{L}$ is an operator of order $2m$ with real-valued basic function coefficients, it is still transversely elliptic because for every $z\in M$, every transverse covector $\xi\in (T^*M\slash \cF^*)\vert_z$ different from $0$,
\begin{equation}\label{eq:4.1}
\sigma(L\circ \overline{L})(z,\xi)=\left|\sigma(L)(z,\xi)\right|^2>0
\end{equation}
where $|\cdot|$ is the module of complex number. For every $z\in M$, the isomorphism between $(T^*M\slash \cF^*)\vert_z$ and $(\nu\cF)^*\vert_z$ permits us to regard $\sigma(L\circ\overline{L})(z,\cdot)$ as an element of $S^{2m,0}_{++}(\nu\cF\vert_z)$. As in each local chart, all coefficients of $L\circ \overline{L}$ are real-valued basic functions, a priori smooth, so $\sigma(L\circ\overline{L})\in C^\infty\left(M,S_{++}^{2m,0}(\nu\cF)\right)$. By Proposition 3.3, $\Phi_{\nu\cF}(\sigma(L\circ\overline{L}))\in C^\infty(M,S_{++}^{0,2}(\nu\cF))$ is a definite-positive $C^\infty(M)$-bilinear form on $C^\infty(M,\nu\cF)$. By remark 1, it rests to prove its holonomy invariance to obtain a transverse metric on $(M,\cF)$. We need only to prove the holonomy invariance in any arbitrary foliated chart.\\
Let $U$ be a foliated chart of coordinates $(x,y)=(x_1,\cdots,x_{n-q},y_1,\cdots,y_q)$ where $L\circ \overline{L}$ is written as 
$$
L\circ \overline{L}=\sum_{|\alpha|\leq 2m} a_{\alpha} \dfrac{\partial^{|\alpha|}}{\partial {y_1}^{\alpha_1} \partial {y_2}^{\alpha_2}\cdots \partial {y_q}^{\alpha_q}}.
$$
On $U$, the quotient map $TM\rightarrow \nu\cF$ renders the local foliated fields ${\partial}/{\partial y_k}\in\frak{X}(U,\cF\vert_U)$ to local basic fields $\left[{\partial}/{\partial y_k}\right]\in l(U,\cF\vert_U)$. Thus, we get a local frame $\mathscr{B}=\left([{\partial}/{\partial y_k}]\right)_{k=1,\cdots,q}$ of $\nu\cF\vert_U$ and the local dual frame $\mathscr{B}^*=\left([dy_k]\right)_{k=1,\cdots,q}$. As $\mathcal{L}_X [{\partial}/{\partial y_k}]=0$ for any $X\in\frak{X}(\cF\vert_U)$, we have $\mathcal{L}_X [d y_k]=0$. Locally on $U$, $\sigma(L\circ\overline{L})$ has the expression
$$
\sigma(L\circ\overline{L})=\sum_{|\alpha|=2m} a_{\alpha} \left[{\partial}/{\partial y_1}\right]^{\alpha_1}\cdots\left[{\partial}/{\partial y_q}\right]^{\alpha_q}
$$
Thus, omitting the variable $z\in U$,
$$
\Psi_{\nu\cF}\Big(\sigma(L\circ\overline{L})\Big)=\sum_{i,j=1}^q \left[\phi\Big(\sigma(L\circ\overline{L})_{\mathscr{B}}\Big)\right]_{ij} [dy_i] [dy_j].
$$
In order to prove the holonomy invariance of $\Psi_{\nu\cF}\Big(\sigma(L\circ\overline{L})\Big)$, it only needs to show each component $\left[\phi(\sigma(L\circ\overline{L})_{\mathscr{B}})\right]_{ij}$ is holonomy invariant in $U$. Under $\mathscr{B}$, let $\xi=[\xi_1,\cdots,\xi_q]^T\in\mathbb{R}^q$, we write
$$
\sigma(L\circ \overline{L})_{\mathscr{B}} (\xi) = \sum_{|\alpha|=2m} a_{\alpha}\xi_1^{\alpha_1}\cdots\xi_q^{\alpha_q}
$$
Since all coefficients $a_\alpha$ are basic functions, each $(i,j)$-component of matrix-valued function $\phi\Big(\sigma(L\circ\overline{L}\Big)_{\mathscr{B}})$ is constructed by basic functions, so $\mathcal{L}_X \left[\phi(\sigma(L\circ\overline{L})_{\mathscr{B}})\right]_{ij}=0$, for any $X\in\frak{X}(\cF\vert_U)$.
\end{proof}

Let us see the special case of $m=1$, i.e. second order transversely elliptic basic differential operator. It is well-known that in the classic theory, we can define the Laplace operator (called also the Laplace-Beltrami operator) acting on space of fucntions on any Riemannian manifold,\cite[Chapter 3]{G}. Or generally the Laplace operator acting on space of sections of vector bundle, \cite{BGV}. Analoguely, in foliated case, if there is a transverse metric $g$ on $(M,\cF)$, we can define the transverse Laplace operator $\Delta_g$ associated to $g$. In each local foliated chart of coordinates $(x,y)$ as above, the transverse Laplace operator has the expression
$$
\Delta_{g} = \sum_{i,j=1}^q \frac{1}{\sqrt{\det g}}\frac{\partial}{\partial y_i}\Big(\sqrt{\det g} g^{ij}\dfrac{\partial}{\partial y_j}\Big).
$$
where by abuse of notation $g$ denotes the metric matrix $[g_{ij}]$ and $[g^{ij}]$ is the inverse matrix of $[g^{ij}]$. We leave readers to check that $\Delta_{g}$ is a real coefficient second order transversely elliptic basic differential operator acting on $C_b(M)$. A basic differential operator acting on $C_b(M)$ is called generalized transverse Laplace operator if it has the same principal part of $\Delta_{g}$ for some transverse metric $g$. We have the following theorem.
\begin{Thm}
Each second order transversely elliptic basic differential operator acting on the space of basic functions with real-valued basic function coefficients is a generalized transverse Laplace operator. 
\end{Thm}

\begin{proof}
Let $L$ be such an operator acting on $C_b(M)$. It suffices to make proof locally. Suppose that in a local foliated chart $U$ as above, $L$ has the expression
$$
L=\sum_{i,j=1}^q a^{ij} \frac{\partial^2}{\partial y_i\partial y_j}+\ \text{lower\ order\ part}.
$$
Locally in $U$, under $\mathscr{B}$, $\mathscr{B}^*$ as above, the principal symbol $\sigma(L)$ has the expression
$$
\sigma(L)=\sum_{i,j=1}^q a^{ij}\left[{\partial}/{\partial y_i}\right]\left[{\partial}/{\partial y_j}\right]
$$
Then, by Proposition \ref{Prop2.5}, $\phi\Big({\sigma(L)}_\mathscr{B}\Big)=[a^{ij}]^{-1}$. So,
$$
\Psi_{\nu\cF}\Big(\sigma(L)\Big)=\sum_{i,j=1}^n a_{ij} \left[d y_i\right]\left[d y_j\right]
$$
where $[a_{ij}]$ is the inverse matrix of $[a^{ij}]$. Locally in $U$, $L$ is a generalized transverse Laplace operator associated to the transverse metric induced by $\Phi_{\nu\cF}(\sigma(L))$ under the pullback $(\nu\cF)^*\rightarrow T^*M$.
\end{proof}

\subsection{Transversely elliptic operator on $\cF$-bundle}
We may generalize Theorem 1 to a transversely elliptic basic differential operator on $\cF$-bundle. Here, we use the notations in \cite[Section2.7]{AKHR}. Let $\mathcal{E}\rightarrow M$ be a complex vector bundle of rank $N$. Recall that $\mathcal{E}$ is said to be a $\cF$-bundle if its associated principal bundle is an $\cF$-bundle, which induces a basic connection $\nabla$ on $\mathcal{E}$. Denote by $C^\infty(\mathcal{E})$ the space of global sections. We say that a section $s\in C^\infty(\mathcal{E})$ is basic if it satisfies the condition $\nabla_X s=0$ for every $X\in \frak{X}(\cF)$. We denote by $\widetilde{\mathcal{E}}$ the sheaf of germs of its basic sections.
\begin{Def}
A basic differential operator of order $m\in\mathbb{N}^*$ on $\mathcal{E}$ is a morphism of sheaves $D:\widetilde{\mathcal{E}}\rightarrow \widetilde{\mathcal{E}}$ such that in a foliated local chart of coordinates $(x,y)=(x_1,\cdots,x_{n-q},y_1,\cdots,y_q)$, $D$ has the expression
$$
D=\sum_{|\alpha|\leq m} a_{\alpha}\dfrac{\partial^{|\alpha|}}{\partial {y_1}^{\alpha_1}\cdots \partial {y_q}^{\alpha_q}}
$$
where $\alpha=(\alpha_1,\cdots,\alpha_q)\in\mathbb{N}^q$, $|\alpha|=\alpha_1+\cdots+\alpha_q$ and $a_{\alpha}$ are $N\times N$ matrices whose components are basic functions. We say that $D$ is transversely elliptic if its principal symbol $\sigma(D)(z,\xi)\in\mathrm{End}(\mathcal{E}\vert_z)$ is an isomorphism for every $z\in M$ and every transverse covector $\xi\in (T^*M\slash \cF^*)\vert_z$ different from $0$.
\end{Def}

\begin{Thm}
Let $M$ be a foliated manifold with a foliation $\cF$ which admits non-constant basic functions. Let $\mathcal{E}\rightarrow M$ be an $\cF$-bundle. If there exists a transversely elliptic basic differential operator on $\mathcal{E}$, then $\cF$ is a Riemannian foliation.
\end{Thm}
\begin{proof}
Let $D$ be such an operator of order $m\in\mathbb{N}^*$. Analoguely to Theorem 1, consider the transversely elliptic operator $D\circ \overline{D}$. As we see locally that all coefficients of $D\circ \overline{D}$ are $N\times N$-matrices of real-valued basic functions, for each $z\in M$ and each transverse covector $\xi\in (T^*M\slash \cF^*)\vert_z$ different from $0$, the determinant of the isomorphism $\sigma(D\circ \overline{D})(z,\xi)$ of $\mathcal{E}\vert_z$ is a non-zero real number. Thus, upto a $\sqrt{-1}$ before $D$ and under the assumption of the connectedness of $M$, we can assume that $\det\Big(\sigma(L\circ \overline{L})(z,\xi)\Big)>0$. So, we have
$$
\det\sigma(L\circ \overline{L}) \in C^{\infty}\left(M,S^{2mN,0}_{++}(\nu\cF)\right)
$$
Followed by the same argument in Theorem 1, we know that $\Psi_{\nu\cF}\Big(\det\sigma(L\circ \overline{L})\Big)$ induces a transverse metric on $(M,\cF)$. Therefore, $\cF$ is a Riemannian foliation.
\end{proof}

\bibliographystyle{plain}

\end{document}